\documentclass[12pt]{article}
\usepackage{latexsym}
 
\newcommand{\proof}{{\noindent \bf Proof. }}
\newtheorem{thm}{Theorem}

\begin{document}
\begin{titlepage}
\title{\bf HAMILTON PATHS WITH LASTING SEPARATION}
{\author{{\bf Emanuela Fachini}
\\{\tt fachini@di.uniroma1.it}
\\''La Sapienza'' University of Rome
\\ ITALY
\and{\bf J\'anos K\"orner}
\thanks{Department of Computer Science, University of Rome, La Sapienza, 
via Salaria 113, 00198 Rome, ITALY}
\\{\tt korner@di.uniroma1.it}
\\''La Sapienza'' University of Rome
\\ ITALY}}

\maketitle
\begin{abstract}

We determine the asymptotics of the largest cardinality of a set of Hamilton paths in the complete graph with vertex set $[n]$ under the condition that for any two of  the paths in the family there is a subpath of length $k$ entirely contained in only one of them and edge--disjoint from the other one.  
\end{abstract}
\end{titlepage}

\section{Introduction}In recent years, in a series of papers we have studied the size of the largest family of Hamilton paths in a fixed complete graph with vertex set $[n]$ such that for any pair of the Hamilton paths from the family their union contains a fixed small subgraph \cite{KMS}, \cite{KMu}, \cite{CFK}. The central question in the last mentioned paper is concerned with the case when the union of the paths must contain 
a cycle of prescribed length. Subsequently, answering a question left open in \cite{CFK}, I. Kov\'acs and D. Solt\'esz \cite{KS} have proved that the size of the largest family of Hamiltonian paths the pairwise union of 
which contains an odd cycle does not decrease if instead of an arbitrary odd cycle the union contains a triangle. (The analogous statement for even cycles vs. cycles of length 4 is false, \cite{CFK}.)  These problems are rooted in zero-error information theory, especially in the graph capacity problem of Claude Shannon \cite{Sh}. The paper \cite{GKV} introduced far-reaching generalisations of Shannon's problem with the aim of answering well-known questions in extremal combinatorics, including R\'enyi's problem of the maximum size of a family of pairwise qualitatively independent partitions of a finite set. The concept of 
capacity of permutations is a further natural extension of graph capacity \cite{KMa} and has led to similar questions about Hamilton paths. 

It has become clear throughout the years that the formula--oriented approach to zero-error capacity is practically hopeless. As a matter of fact, even the question about the exact value of the zero-error capacity of odd cycles of length larger than 5 and of their complements is so far from being solved that not even a conjectured value of these capacities has been put forward. For this reason it seems more promising to ask qualitative questions about graph capacities. It was explained in \cite{KMaS} that the question about the capacity of infinite graphs generalises Shannon's original question, more precisely, it explicitly contains the more technical problem of graph capacity for fixed composition codes. To fix ideas, we will suppose that the set of vertices of our graphs is the set of natural numbers. It is clear and simple to see that if an infinite graph has a finite chromatic number, then it has a finite capacity. One can consider the graph induced by the {\it co-normal} power graph $G^n$ on the permutations of $[n]$, seen as $n$--length sequences of vertices of $G.$ 
If the chromatic number of $G$ is infinite, it is very challenging to investigate how fast the clique number of $G$ increases with $n$, more precisely, we ask whether the increase is only exponential or larger. 
To encompass a very large set of these problems, we can think of permutations as directed Hamilton paths with vertex set $[n]$. We ask similar questions for non-directed Hamilton paths as well. The questions about graph capacity are formulated in terms of pairwise difference in particular pairs of natural numbers. Considering permutations in their graph representation by Hamilton paths, other, less local conditions of difference become natural. 
The questions studied in the present paper are intimately connected to those raised in \cite{CFK}, even though the conditions studied here are not in terms of subgraphs in the pairwise union of Hamilton paths. Our main concern is to understand the rough order of magnitude of the cardinality of the largest set of paths satisfying various conditions in order to develop an intuition for their hierarchy. 

Somewhat surprisingly, our present conditions allow for large sets of superexponential size. 

\section{A simple problem}

Let us start with the case $k=2$ since already this exhibits the main features of what we see in general. Let $M(n,2)$ be the largest cardinality of a family of Hamilton paths in the complete graph $K_n$ with vertex set $[n]$ such that for any two of them there is a path of 3 vertices with both edges contained in the same path and missing in the other one. We claim

\begin{thm}\label{thm:two}

$${{\lfloor {n \over 2} \rfloor}! \over  2^{\lfloor  n/2 \rfloor}}    \leq M(n,2) \leq 2^{\lceil n/2 \rceil}{\lceil {n \over 2}\rceil}!$$
\end{thm}

\proof

The lower bound is an immediate consequence of a result in \cite{KM}. In \cite{KM}, slightly improving on a previous result from \cite{KMu}, it was shown that there exists a family of Hamilton paths in 
$K_n$ with the property that the union of any two of the paths from the family has a vertex of degree 4. Now, it is immediate that such a family has the property we require in this Theorem. In fact, if a vertex has degree 4 in the union of two paths, then the two paths have no common edge incident to this vertex, whence the two edges coming from the same path are consecutive in that path and do not belong to the other path.
Hence the cardinality of the family constructed in \cite{KM} is a lower bound to $M(n,2).$

For the upper bound, consider first the case of $n$ even. Fix a perfect matching in $K_n$ and consider all the Hamilton paths containing all its edges. These various paths are defined by specifying an arbitrary order in which the edges of the matching appear in the Hamilton path alongside with the order of appearance of the two vertices for any of these edges. Then, clearly, we obtain $2^{n/2-1}{n \over 2}!$ Hamilton path in this manner and these are exactly those that contain the given matching. Every Hamilton path in $K_n$ contains exactly one perfect matching. Clearly, two paths containing the same perfect matching cannot have a ''private" path on 3 vertices, since every second edge of a Hamilton path is an edge from the fixed matching, and therefore is a common edge of both of the paths. This implies that an optimal construction of $M(n,2)$ Hamilton paths cannot have more paths than there are perfect matchings in $K_n$ and this gives our upper bound in case of even $n.$ The case of odd $n$ can be treated similarly. For the upper bound, we have to note that even though every Hamilton path defines two edge--disjoint near--perfect matchings, the number of these near--perfect matchings is at most
$${n! \over  2^{\lfloor n/2 \rfloor} {\lfloor n/2 \rfloor}!} \leq 2^{\lceil n/2 \rceil}{\lceil {n \over 2}\rceil}!$$
and it is still true that two paths containing the same near--perfect matching cannot satisfy our condition, 
which gives the upper bound.

\hfill$\Box$

We can adapt our previous reasoning to the case of $k>2.$ Let $M(n,k)$ be the maximum size of a family of Hamilton paths in $K_n$ such that for any pair of the Hamilton paths in the family there is a subpath of 
$k$ edges in one of the paths that is edge--disjoint from the whole of the other path. (This subpath is obviously induced by $k+1$ vertices.) We treat this problem in two separate results according to the parity of $k.$

\begin{thm}\label{thm: fo}

If $k>2$ is even and $n$ is a multiple of $k$, we have 

$$(n/k)! \leq M(n,k) \leq 3^n(n/k)!$$

\end{thm}

\proof

In order to prove the lower bound, we consider the complete bipartite graph $K_{n/2.n/2}$ with its independent sets $A$ and $B$ of equal size. We fix an order of the vertices in $A$ while in $B$ we first partition the vertices into disjoint ordered $(k-1)$-tuples. During the construction process these tuples and their internal order will stay fixed. On the other hand we will vary the relative order of the different tuples in every possible way. For each of the $(n/k)!$ of the orders of the elements of $B$ so obtained, we define a Hamilton path in the whole bipartite graph as follows. Let all our Hamilton paths start with the first vertex 
from $A.$ From here the path goes to the first vertex of the first $(k-1)$-tuple from $B.$ From the latter the path continues to the second vertex of $A,$ and then goes back to the second vertex of the first $(k-1)$-tuple of $B.$ In other words, after having chosen our permutation of the vertices of $B,$ the construction is defined precisely as in the previous proof of Theorem \ref{thm:two}. The number of the paths so obtained 
equals the number of the considered permutations of $B$ which is $(n/k)!.$ Since the relative order of the$(k-1)$-tuples is different, therefore, wherever the two permutations of $B$ differ, both of them will generate a path of k edges, all of which incident to one of the vertices of the corresponding $(k-1)$-tuple from $B.$  Furthermore, since the two tuples are disjoint, the two paths of $k$ edges do not have common edges, and none of these edges appears elsewhere in the union of the two paths.

In order to establish our upper bound, consider the set of vertex-disjoint edges $$\{ki+1,ki+2\}\quad  \hbox{for}\;  i=0,1,2,\dots, n/k-1$$  of the complete graph $K_n.$ To any set of edges such as this we associate all the Hamilton paths that visit these edges in any order 
and with the further restriction that between any pair of successively visited edges the path traverses exactly $k-1$ arbitrary new vertices not belonging to any of the edges. It is then clear that no pair of Hamilton paths associated to a fixed edge set in this manner satisfies our pairwise condition since each of them has every $k$'th of its edges in common with all of  the other ones. The number of such paths is clearly 
$$(n/k)! 2^{n/k}(n-{2n \over k})!$$ 
and of all these at most one can belong to our set of Hamilton paths satisfying the pairwise condition. Hence the total number of paths in our family is at most
$${n! \over 2(n/k)! 2^{n/k}(n-{2n \over k})!}\leq $$
$$\leq {n!\over (n/k)!2^{n/k}(n-{2n \over k})!}= {n!(n/k)! \over  (n/k)!2^{n/k}(n-{2n \over k})! (n/k)!} < 3^n(n/k)!2^{-n/k}< 3^n(n/k)!$$

\hfill$\Box$

If $k$ is odd, things change a little bit in the construction. 

\begin{thm}\label{thm: thre}

If $k$ is odd and $n$ is a multiple of $k$, 
$$(n/k-1)! \leq M(n,k) \leq 3^n(n/k)!$$
\end{thm}

\proof
In order to prove the lower bound, set $a={\lfloor k/2 \rfloor \over k}n$ and consider the bipartite complete graph $K_{a, n-a}$ with vertex set $[n].$ Let the two independent sets be $A$ and $B$ with $|A|=a.$ For the ret of the construction, let us fix an order of the elements of $A.$ 
Further, let us partition these elements into disjoint $\lfloor k/2 \rfloor$-tuples of consecutive elements. In the set $B$ let us partition the elements in an arbitrary manner into disjoint $\lceil k/2 \rceil$-tuples and let us fix an arbitrary order of the $\lceil k/2 \rceil$ elements within each of these groups. Let us now permute the $\lceil k/2 \rceil$-tuples of $B$ arbitrarily except for one that we will keep fixed as the first one. Each of these permutations will give rise to a different Hamilton path in our bipartite complete graph as follows. Each of these paths will start with the first element of the first $\lceil k/2 \rceil)$-tuple of $B$ and will alternate between the elements of the first $\lceil k/2 \rceil$-tuple of $B$ and the first  $\lfloor k/2 \rfloor$-tuple of $A.$ 
Once the path arrives to the last element of the first   $\lceil k/2 \rceil$-tuple of $B$ it interrupts alternation to continue with the first vertex of the second  $\lceil k/2 \rceil$-tuple of $B$. From here the path goes to the first element of the second $\lfloor k/2 \rfloor$-tuple of $A$ to resume alternation. We keep repeating this procedure until all the vertices are covered. The number of the corresponding Hamilton paths is clearly 
the factorial of
$${|B| \over   \lceil k/2 \rceil}={n\lceil k/2 \rceil \over k} ( \lceil k/2 \rceil)^{-1} ={n \over k}$$
This factorial gives the lower bound, once we realise that any pair of the Hamilton paths so obtained satisfies our pairwise condition. To see that any pair of these paths satisfy our pairwise condition, let us look at the corresponding permutations of the elements of $B.$ We will refer to these as the first and the second permutation.
Set $t= \lceil k/2 \rceil$ and let $(a_1, a_2, \dots a_t)$ be the first $t$-tuple of elements in the first permutation that differs, and in fact, is disjoint from the $t$-tuple in the same positions in the second permutation. Let further $x$ be the first vertex of the $t$-tuple following  $(b_1, b_2, \dots b_t)$ in the first permutation and let $(y_1, y_2, \dots y_{k-t}$ be the corresponding $k-t$-tuple in $A.$ Then we see that the consecutive edges $( (b_1, y_1), (y_1, b_2), \dots, (y_{k-t},b_t)$,  and $(b_t,x)$ define a subpath 
in the Hamilton path corresponding to the first permutation and none of these edges belongs to the Hamilton path corresponding to the second permutation.

The proof of the upper bound is more straightforward. We chose the set of edges $E=\{\{ki+1,ki+2\}, \quad \hbox{for}\;  i=0,1,2, \dots, n/k-1\}.$ We consider the set of all the Hamilton paths that visit these edges in any order under the restriction that between any two edges the path passes exactly $k-2$ different vertices not contained in any of them. Each of these vertices is visited just once. Since every $k$'th edge in each such path is from our set $E$ and thus belongs to all of the paths, it follows that no pair of paths from this set satisfies our condition. On the other hand, the number of paths corresponding to a fixed set $E$ in this manner is
$$2^{n/k}\cdot (n/k)!(n-{2n\over k})!.$$
This leads to the  required upper bound in the same way as in the previous proof.

\hfill$\Box$

\section{Kernel and product structure}

Interestingly enough, our near--optimal constructions in this paper are in terms of kernel structures. In other words, the family of Hamilton paths we thus suspect to be almost optimal in various problems considered has a so--called {\em kernel structure.} This means that all the paths in the construction have a fixed projection while the rest varies arbitrarily on a vertex set of size linear in $n.$ This is the reason why these constructions are superexponential in size. As a matter of fact, we fix a linearly ordered set of $n/2$ vertices and consider all the Hamilton paths in $K_{n}$ which have these vertices as odd--indexed vertices in the order of transition through the vertices of 
$K_n.$ All our constructions have a similar structure. In a sense, the structures we build can be considered as a Cartesian product. More generally, it is interesting to analyse when it is that in a problem in 
extremal combinatorics, the extremal solutions have a Cartesian product structure. In Shannon's graph capacity problem, in all the solved basic cases, the optimal construction has a product structure. (This is so even for the famous pentagon graph, as it was proved by Lov\'asz \cite{L}.) 

On the other hand, in most of the analogous problems about families of Hamilton paths product structures will not give any non-trivial construction \cite{KM}, \cite{CFK}. 

\section{Open problems}

We might wonder about the size $F(n,k)$ of the largest family of subgraphs of $K_n$ such that, for any two of them, one contains a not necessarily induced path of length $k$ that is completely edge-disjoint from the other 
graph. Is it true that the optimum (or the near-optimum) is reached for families of Hamilton paths?

Questions about the difference or the symmetric difference of Hamilton paths have not been considered before in the present information-theoretic framework. We suggest to look at them more closely. In this spirit let us examine 
and denote by $L(n,k)$ the largest cardinality of a family of Hamilton paths in $K_n$ such that for any pair of its members there is a $k$-matching, i. e., a set of $k$ vertex-disjoint edges contained in just one of 
the two paths and edge-disjoint from the other one. As before, for every fixed $k$ we would like to determine the asymptotics of $L(n,k)$ in $n.$ Is this asymptotics the same as for $M(n,k)$? To conclude, let us look at triangles. Let us denote by $D(n,3)$ the largest cardinality of a family of subgraphs of $K_n$ such that for any two of its member graphs there is a triangle contained in exactly one of the two graphs and edge-disjoint from the other one. If we restrict attention to induced subgraphs of our complete graph and represent them by the $n$-length binary vectors which are their characteristic vectors, a large family of these can be obtained from error-correcting codes. Their existence is guaranteed by the Gilbert-Varshamov bound \cite{Gi}, \cite{V}, giving, for $n$ large enough, the lower bound 
$$D(n,3)\geq {2^n \over n^6}$$
We have no non-trivial upper bound. 
\section{Acknowledgement} 

We are grateful to an anonymous referee for having spotted an error in the original proof of Theorem \ref{thm:two} and various useful suggestions for rewriting the paper. 

\end{document}